\documentclass[11pt,reqno]{amsart}

\usepackage[T1]{fontenc}
\usepackage{lmodern}
\usepackage{microtype}
\usepackage{amsmath,amssymb,amsthm,mathtools}
\usepackage[margin=1.02in]{geometry}
\usepackage{xcolor}
\usepackage[colorlinks=true,linkcolor=blue!60!black,citecolor=blue!60!black,urlcolor=blue!60!black]{hyperref}

\newtheorem{theorem}{Theorem}[section]
\newtheorem{corollary}[theorem]{Corollary}

\newcommand{\C}{\mathbb{C}}

\title[Bourin's Crossed-Heinz Inequality in the Operator and Frobenius Norms]{Bourin's Crossed-Heinz Inequality in the Operator and Frobenius Norms}

\author[Yulong Tian]{Yulong Tian$^{*}$}
\address{Department of Mathematics, Northeast Forestry University, Harbin, P. R. China}
\email{tylmath@nefu.edu.cn}
\thanks{$^{*}$Corresponding author.}

\author[Jinli Xu]{Jinli Xu}
\address{Department of Mathematics, Northeast Forestry University, Harbin, P. R. China}
\email{jclixv@qq.com}


\begin{document}

\begin{abstract}
For positive semidefinite matrices $A, B \in M_n(\mathbb{C})$ and $0 < t < 1$, we establish Bourin's crossed-Heinz inequality for every dimension $n \ge 1$ in both the operator norm and the Frobenius norm. In the operator norm, we prove the sharp inequality
\[
\left\|A^t B^{1-t} + B^t A^{1-t}\right\|_{\mathrm{op}} \le \left\|A + B\right\|_{\mathrm{op}}
\]
by normalizing $\|A + B\|_{\mathrm{op}}$ to $2$ and decomposing the normalized crossed-Heinz sum into $P_t P_{1-t} + Q$, where $P_t$, $P_{1-t}$, and $Q$ are self-adjoint contractions. In the Frobenius norm, we establish the sharp quantitative inequality
\[
\left\|A^t B^{1-t} + B^t A^{1-t}\right\|_{\mathrm{F}}^{\,2} + 4t(1-t)\left\|A - B\right\|_{\mathrm{F}}^{\,2} \le \left\|A + B\right\|_{\mathrm{F}}^{\,2}
\]
via an exact nonnegative remainder decomposition, showing in particular that $\|A^t B^{1-t} + B^t A^{1-t}\|_{\mathrm{F}} \le \|A + B\|_{\mathrm{F}}$ holds with equality if and only if $A = B$.
\end{abstract}

\maketitle

\section{Introduction and Main Results}

Let $M_n(\C)$ denote the space of $n \times n$ complex matrices, with the identity matrix $I$. We write $\Vert\cdot\Vert_{\mathrm{op}}$ and $\Vert\cdot\Vert_F$ for the operator and Frobenius norms on $M_n(\mathbb{C})$, respectively. For positive semidefinite $A, B \in M_n(\mathbb{C})$ and $t \in (0, 1)$, the classical Heinz inequality of Bhatia and Davis~\cite{BD93} (see also~\cite[Chapter~IX]{Bha97}) states that $H_t(A, B) := A^t B^{1-t} + A^{1-t} B^t$ satisfies
\[
\bigl|\!\bigl|\!\bigl| H_t(A, B) \bigr|\!\bigr|\!\bigr| \le \bigl|\!\bigl|\!\bigl| A + B \bigr|\!\bigr|\!\bigr|
\]
for every unitarily invariant norm $\lvert\!\lvert\!\lvert \cdot \rvert\!\rvert\!\rvert$ on $M_n(\mathbb{C})$. In~\cite[Question~4]{Bou09}, Bourin asked whether the corresponding bound remains valid for the crossed-Heinz sum $C_t(A, B) := A^t B^{1-t} + B^t A^{1-t}$, with the factors in the second term reversed:
\begin{equation}\label{eq:bourin}
\bigl|\!\bigl|\!\bigl| C_t(A, B) \bigr|\!\bigr|\!\bigr| \le \bigl|\!\bigl|\!\bigl| A + B \bigr|\!\bigr|\!\bigr|.
\end{equation}

A natural route toward~\eqref{eq:bourin} would be to dominate $C_t(A, B)$ by $H_t(A, B)$, yet this comparison fails in general in both the operator and Frobenius norms. In the operator norm, Bottazzi, Elencwajg, Larotonda, and Varela~\cite[Example~3.1]{BELV15} constructed $3 \times 3$ positive definite matrices with $\|C_t(A, B)\|_{\mathrm{op}} > \|H_t(A, B)\|_{\mathrm{op}}$. In the Frobenius norm, Bhatia~\cite{Bha14} established the stronger bound $\|C_t(A, B)\|_{\mathrm{F}} \le \|H_t(A, B)\|_{\mathrm{F}}$ for $t \in [\tfrac{1}{4}, \tfrac{3}{4}]$, whereas counterexamples by Plevnik~\cite{Ple16} and later by Carlen and Lieb~\cite{CL22} show that this inequality does not extend to the full interval $(0, 1)$.

Following a series of partial and quantitative bounds by Hayajneh, Hayajneh, and Kittaneh~\cite{HHK17,HHK21,HHK22,HHK23a,HHK23b}, Kittaneh and Ricard~\cite{KR25} established~\eqref{eq:bourin} for all unitarily invariant norms and all dimensions when $t \in [\tfrac{1}{4}, \tfrac{3}{4}]$. Most recently, Sababheh and Bourin~\cite[Theorem~3.6]{SB26} proved $\lvert\!\lvert\!\lvert C_t(A, B) \rvert\!\rvert\!\rvert \le \lvert\!\lvert\!\lvert H_t(A, B) \rvert\!\rvert\!\rvert \le \lvert\!\lvert\!\lvert A + B \rvert\!\rvert\!\rvert$ across the full interval $0 < t < 1$ for $2 \times 2$ matrices, while noting that for $n \ge 3$ and $t \notin [\tfrac{1}{4}, \tfrac{3}{4}]$ Bourin's problem remains open.

We establish Bourin's inequality~\eqref{eq:bourin} in both the operator norm and the Frobenius norm for all $n \ge 1$ and $0 < t < 1$.

\begin{theorem}\label{thm:main-op}
Let $n \ge 1$ be an integer, let $A, B \in M_n(\mathbb{C})$ be positive semidefinite matrices, and let $0 < t < 1$. Then
\begin{equation}\label{eq:main-op}
\left\|A^t B^{1-t} + B^t A^{1-t}\right\|_{\mathrm{op}} \le \left\|A + B\right\|_{\mathrm{op}},
\end{equation}
where the constant $1$ on the right-hand side is sharp.
\end{theorem}

\begin{theorem}\label{thm:main-frob}
Let $n \ge 1$ be an integer, let $A, B \in M_n(\mathbb{C})$ be positive semidefinite matrices, and let $0 < t < 1$. Then
\begin{equation}\label{eq:quant-ineq}
\left\|A^t B^{1-t} + B^t A^{1-t}\right\|_{\mathrm{F}}^{\,2} + 4t(1-t)\left\|A - B\right\|_{\mathrm{F}}^{\,2} \le \left\|A + B\right\|_{\mathrm{F}}^{\,2},
\end{equation}
where for each fixed $t \in (0, 1)$ the constant $4t(1-t)$ is sharp. In particular,
\begin{equation}\label{eq:bourin-frobenius}
\left\|A^t B^{1-t} + B^t A^{1-t}\right\|_{\mathrm{F}} \le \left\|A + B\right\|_{\mathrm{F}},
\end{equation}
and for every fixed $t \in (0, 1)$, equality holds in~\eqref{eq:bourin-frobenius} if and only if $A = B$.
\end{theorem}

\section{Contraction Decomposition for the Operator Norm}\label{sec:op}

For Hermitian matrices $A, B \in M_n(\mathbb{C})$, we write $A \le B$ (or $B \ge A$) if $B - A$ is positive semidefinite. Our strategy is to normalize $\|A + B\|_{\mathrm{op}}$ to $2$ and express the normalized crossed-Heinz sum as the sum of a product of two self-adjoint contractions and a third self-adjoint contraction.

For $0 < t < 1$, define
\begin{equation}\label{eq:d-def}
d(x) := (1 - x^t)(1 - x^{1-t}) = 1 - x^t - x^{1-t} + x \qquad (x \ge 0).
\end{equation}
For positive semidefinite matrices $A, B \in M_n(\mathbb{C})$, set
\begin{equation}\label{eq:PQ-def}
P_r := A^r + B^r - I \quad (r \in \{t, 1-t\}), \qquad Q := I - d(A) - d(B).
\end{equation}

\begin{theorem}\label{thm:decomp}
Let $A, B \in M_n(\mathbb{C})$ be positive semidefinite matrices satisfying $A + B \le 2I$, and let $0 < t < 1$. With the notation above, one has $0 \le d(x) \le 1$ for all $x \in [0, 2]$,
\begin{equation}\label{eq:loewner-bounds}
-I \le P_t \le I, \qquad -I \le P_{1-t} \le I, \qquad -I \le Q \le I,
\end{equation}
and
\begin{equation}\label{eq:decomp}
C_t(A, B) = P_t P_{1-t} + Q.
\end{equation}
Consequently, $P_t$, $P_{1-t}$, and $Q$ are self-adjoint contractions, and $\|C_t(A, B)\|_{\mathrm{op}} \le 2$.
\end{theorem}

\begin{proof}[Proof of Theorems~\ref{thm:decomp} and~\ref{thm:main-op}]
For $r \in \{t, 1-t\}$, the scalar function $x \mapsto x^r$ is concave on $[0, \infty)$, so its tangent line at $1$ gives $x^r \le (1 - r) + rx$ for all $x \ge 0$. Hence $A^r \le (1 - r)I + rA$ and $B^r \le (1 - r)I + rB$, and adding the two bounds yields
\[
0 \le A^r + B^r \le 2(1 - r)I + r(A + B) \le 2I.
\]
Thus $-I \le P_r \le I$. Since $P_r$ is self-adjoint, $\|P_r\|_{\mathrm{op}} \le 1$.

Next, we verify that $0 \le d(x) \le 1$ on $[0, 2]$. For $0 \le x \le 1$, both $1 - x^t$ and $1 - x^{1-t}$ lie in $[0, 1]$. For $1 \le x \le 2$, we have $1 \le x^t, x^{1-t} \le 2$, so $d(x) = (x^t - 1)(x^{1-t} - 1)$ is again a product of two numbers in $[0, 1]$. Since $0 \le A \le A + B \le 2I$ and $0 \le B \le A + B \le 2I$, the spectra of $A$ and $B$ lie in $[0, 2]$, whence $0 \le d(A) \le I$ and $0 \le d(B) \le I$. Adding gives $0 \le d(A) + d(B) \le 2I$, so $Q = I - d(A) - d(B)$ satisfies $-I \le Q \le I$. Thus $Q$ is a self-adjoint contraction, establishing~\eqref{eq:loewner-bounds}.

Expanding the product gives
\[
(A^t + B^t)(A^{1-t} + B^{1-t}) = A^t A^{1-t} + A^t B^{1-t} + B^t A^{1-t} + B^t B^{1-t} = A + B + C_t(A, B).
\]
By~\eqref{eq:d-def} and~\eqref{eq:PQ-def}, we obtain
\[
\begin{aligned}
P_t P_{1-t} &= A + B + C_t(A, B) - (A^t + B^t + A^{1-t} + B^{1-t}) + I, \\
Q &= -I + (A^t + B^t + A^{1-t} + B^{1-t}) - (A + B).
\end{aligned}
\]
Adding these two identities yields~\eqref{eq:decomp}. By the triangle inequality and submultiplicativity of the operator norm,
\[
\|C_t(A, B)\|_{\mathrm{op}} \le \|P_t\|_{\mathrm{op}}\,\|P_{1-t}\|_{\mathrm{op}} + \|Q\|_{\mathrm{op}} \le 1 \cdot 1 + 1 = 2,
\]
which proves Theorem~\ref{thm:decomp}.

For Theorem~\ref{thm:main-op}, if $A + B = 0$, then $A = B = 0$ and~\eqref{eq:main-op} is trivial. Otherwise, set
\[
\lambda := \frac{\|A + B\|_{\mathrm{op}}}{2} > 0, \qquad \widehat{A} := \frac{A}{\lambda}, \qquad \widehat{B} := \frac{B}{\lambda}.
\]
Then $\widehat{A}, \widehat{B} \ge 0$ and $\|\widehat{A} + \widehat{B}\|_{\mathrm{op}} = 2$, so $\widehat{A} + \widehat{B} \le 2I$. Since $C_t(A, B) = \lambda\,C_t(\widehat{A}, \widehat{B})$, Theorem~\ref{thm:decomp} gives
\[
\|C_t(A, B)\|_{\mathrm{op}} = \lambda\,\|C_t(\widehat{A}, \widehat{B})\|_{\mathrm{op}} \le 2\lambda = \|A + B\|_{\mathrm{op}},
\]
proving~\eqref{eq:main-op}. Finally, taking $A = B \neq 0$ gives $C_t(A, A) = 2A = A + B$, showing that the constant $1$ in~\eqref{eq:main-op} is sharp.
\end{proof}

Since the proof of Theorem~\ref{thm:decomp} relies only on continuous functional calculus and basic order and norm properties, replacing $I$ by the unit extends the argument verbatim to arbitrary unital $C^*$-algebras.

\begin{corollary}\label{cor:infinite-dim}
Let $\mathcal{A}$ be a unital $C^*$-algebra with positive cone $\mathcal{A}_+$, and let $0 < t < 1$. Then
\[
\left\|A^t B^{1-t} + B^t A^{1-t}\right\| \le \left\|A + B\right\|, \qquad A, B \in \mathcal{A}_+.
\]
In particular, for bounded positive operators on any complex Hilbert space $\mathcal{H}$,
\[
\left\|A^t B^{1-t} + B^t A^{1-t}\right\| \le \left\|A + B\right\|, \qquad A, B \in B(\mathcal{H})_+.
\]
\end{corollary}

\section{Exact Remainder Decomposition for the Frobenius Norm}\label{sec:decomposition}

For $0 < t < 1$ and positive semidefinite matrices $A, B \in M_n(\mathbb{C})$, define the remainder matrix $R_t(A, B) \in M_n(\mathbb{C})$ by
\[
R_t(A, B) := (2t - 1)(A - B) - A^t B^{1-t} + B^t A^{1-t}.
\]
Choose arbitrary orthonormal eigenbases $\{e_1, \dots, e_n\}$ and $\{f_1, \dots, f_n\}$ of $\mathbb{C}^n$ for $A$ and $B$, respectively, with
\[
A e_i = a_i e_i, \qquad B f_j = b_j f_j, \qquad a_i, b_j \ge 0 \qquad (1 \le i, j \le n).
\]
For nonnegative real scalars $a, b \ge 0$, define
\begin{align*}
k_t(a, b) &:= a^t b^{1-t} - a^{1-t} b^t, \\
h_t(a, b) &:= k_t(a, b)\bigl[(2t - 1)(a - b) - k_t(a, b)\bigr].
\end{align*}

To prove Theorem~\ref{thm:main-frob}, we express the gap in~\eqref{eq:quant-ineq} as a sum of nonnegative terms.

\begin{theorem}\label{thm:identity}
Let $A, B \in M_n(\mathbb{C})$ be positive semidefinite matrices and let $0 < t < 1$. With the notation above, one has $h_t(a, b) \ge 0$ for all $a, b \ge 0$, and
\begin{equation}\label{eq:exact-decomp}
\|A + B\|_{\mathrm{F}}^{\,2} - \|C_t(A, B)\|_{\mathrm{F}}^{\,2}
= 4t(1-t)\|A - B\|_{\mathrm{F}}^{\,2} + 2\sum_{i,j=1}^{n} |e_i^* f_j|^2\, h_t(a_i, b_j) + \|R_t(A, B)\|_{\mathrm{F}}^{\,2}.
\end{equation}
\end{theorem}

\begin{proof}[Proof of Theorems~\ref{thm:identity} and~\ref{thm:main-frob}]
We first verify that $h_t(a, b) \ge 0$ for all $a, b \ge 0$. By symmetry, we may assume $a \ge b \ge 0$ and $t \ge \tfrac{1}{2}$. If $b = 0$, then $h_t(a, 0) = 0$. If $b > 0$, set $z := \tfrac{1}{2}\log(a/b) \ge 0$, so that
\[
k_t(a, b) = 2\sqrt{ab}\,\sinh\bigl((2t-1)z\bigr), \qquad a - b = 2\sqrt{ab}\,\sinh z.
\]
By convexity of $\sinh$ on $[0, \infty)$ and $\sinh 0 = 0$, we have $0 \le \sinh\bigl((2t-1)z\bigr) \le (2t-1)\sinh z$ for $0 \le 2t - 1 < 1$ and $z \ge 0$. Hence $0 \le k_t(a, b) \le (2t-1)(a - b)$, and therefore $h_t(a, b) \ge 0$.

To establish~\eqref{eq:exact-decomp}, introduce the auxiliary matrices
\[
U := A^t B^{1-t} - \frac{2t-1}{2}(A - B), \qquad V := A^{1-t} B^t + \frac{2t-1}{2}(A - B).
\]
Then
\[
U + V^* = C_t(A, B), \qquad V^* - U = R_t(A, B).
\]
Applying the parallelogram identity for the Frobenius norm together with $\|V^*\|_{\mathrm{F}} = \|V\|_{\mathrm{F}}$ yields
\begin{equation}\label{eq:parallelogram}
\|C_t(A, B)\|_{\mathrm{F}}^{\,2} + \|R_t(A, B)\|_{\mathrm{F}}^{\,2} = 2\|U\|_{\mathrm{F}}^{\,2} + 2\|V\|_{\mathrm{F}}^{\,2}.
\end{equation}

Next, with respect to the orthonormal eigenbases $\{e_i\}$ and $\{f_j\}$, every matrix $Z \in M_n(\mathbb{C})$ satisfies
\begin{equation}\label{eq:mixed-norm}
\|Z\|_{\mathrm{F}}^{\,2} = \sum_{i,j=1}^{n} |e_i^* Z f_j|^2.
\end{equation}
In these bases,
\[
\begin{alignedat}{2}
e_i^* U f_j &= \Bigl(a_i^t b_j^{1-t} - \frac{2t-1}{2}(a_i - b_j)\Bigr) e_i^* f_j,
\qquad &e_i^* V f_j &= \Bigl(a_i^{1-t} b_j^t + \frac{2t-1}{2}(a_i - b_j)\Bigr) e_i^* f_j, \\[5pt]
e_i^* (A + B) f_j &= (a_i + b_j) e_i^* f_j,
&e_i^* (A - B) f_j &= (a_i - b_j) e_i^* f_j.
\end{alignedat}
\]

For any pair of nonnegative scalars $a, b \ge 0$, write $k := k_t(a, b)$. A direct calculation gives
\begin{align*}
&(a + b)^2 - 2\Bigl(a^t b^{1-t} - \frac{2t-1}{2}(a - b)\Bigr)^2 - 2\Bigl(a^{1-t} b^t + \frac{2t-1}{2}(a - b)\Bigr)^2 \\
&\qquad = (a + b)^2 - (a^t b^{1-t} + a^{1-t} b^t)^2 - \bigl[k - (2t-1)(a - b)\bigr]^2 \\
&\qquad = (a - b)^2 - k^2 - \bigl[k - (2t-1)(a - b)\bigr]^2 \\
&\qquad = \bigl[1 - (2t-1)^2\bigr](a - b)^2 + 2k\bigl[(2t-1)(a - b) - k\bigr] \\
&\qquad = 4t(1-t)(a - b)^2 + 2h_t(a, b).
\end{align*}
Applying this identity to $(a_i, b_j)$, multiplying by $|e_i^* f_j|^2$, summing over $1 \le i, j \le n$, and using~\eqref{eq:parallelogram}--\eqref{eq:mixed-norm} gives~\eqref{eq:exact-decomp}.

Since $|e_i^* f_j|^2, h_t(a_i, b_j) \ge 0$,~\eqref{eq:exact-decomp} gives~\eqref{eq:quant-ineq}, and hence~\eqref{eq:bourin-frobenius}. Since $4t(1-t) > 0$, equality in~\eqref{eq:bourin-frobenius} holds if and only if $A = B$.

Finally, for $A = aI$ and $B = bI$ with distinct scalars $a, b > 0$, we have
\begin{align*}
\frac{\|A + B\|_{\mathrm{F}}^{\,2} - \|C_t(A, B)\|_{\mathrm{F}}^{\,2}}{\|A - B\|_{\mathrm{F}}^{\,2}}
&= 1 - \left(\frac{a^t b^{1-t} - a^{1-t} b^t}{a - b}\right)^2 \\
&\longrightarrow 1 - (2t - 1)^2 = 4t(1-t) \qquad (a \to b > 0),
\end{align*}
showing that the constant $4t(1-t)$ in~\eqref{eq:quant-ineq} is sharp.
\end{proof}

\end{document}